\documentclass[12pt]{article}
\usepackage{amsmath, amssymb, amsthm, amsfonts}
\usepackage{amssymb}
\usepackage{amscd}
\usepackage{verbatim}
\begin{document}
\newcommand{\dx}{\,\mathrm{d}x}
\newcommand{\dy}{\,\mathrm{d}y}
\newcommand{\dz}{\,\mathrm{d}z}
\newcommand{\dt}{\,\mathrm{d}t}
\newcommand{\core}{C_0^{\infty}(\Omega)}
\newcommand{\sob}{W^{1,p}(\Omega)}
\newcommand{\sobloc}{W^{1,p}_{\mathrm{loc}}(\Omega)}
\newcommand{\merhav}{{\mathcal D}^{1,p}}
\newcommand{\be}{\begin{equation}}
\newcommand{\ee}{\end{equation}}
\newcommand{\mysection}[1]{\section{#1}\setcounter{equation}{0}}
%%%%%%%%%%%%%%%
\newlength{\wex}  \newlength{\hex}
\newcommand{\understack}[3]{%
 \settowidth{\wex}{\mbox{$#3$}} \settoheight{\hex}{\mbox{$#1$}}
 \hspace{\wex}
 \raisebox{-1.2\hex}{\makebox[-\wex][c]{$#2$}}
 \makebox[\wex][c]{$#1$}
  }%
%%%%%%%%%%%%%
\newcommand{\bea}{\begin{eqnarray}}
\newcommand{\eea}{\end{eqnarray}}
\newcommand{\bean}{\begin{eqnarray*}}
\newcommand{\eean}{\end{eqnarray*}}
\newcommand{\thkl}{\rule[-.5mm]{.3mm}{3mm}}
%%%%%%%%%%%%%%%%%%%%%%%%%%%
\newcommand{\cw}{\stackrel{\rightharpoonup}{\rightharpoonup}}
\newcommand{\id}{\operatorname{id}}
\newcommand{\supp}{\operatorname{supp}}
\newcommand{\wlim}{\mbox{ w-lim }}
\newcommand{\mymu}{{x_N^{-p_*}}}
\newcommand{\Reals}{{\mathbb R}}
\newcommand{\N}{{\mathbb N}}
\newcommand{\Z}{{\mathbb Z}}
\newcommand{\Q}{{\mathbb Q}}
\newcommand{\CC}{\mathbb{C}}
\newcommand{\Hess}{\mathrm{Hess}}
%%%%%%%%%%%%%%%%%
\newcommand{\Real}{{\mathbb R}}
\newcommand{\R}{{\mathbb R}}
\newcommand{\Rn}{\mathbb{R}^n}
\newcommand{\Zn}{\Z^n}
\newcommand{\Nat}{{\mathbb N}}
%%%%%%%%%%%%%%%%%%%%%%
\newcommand{\abs}[1]{\lvert#1\rvert}
\newcommand{\Green}[4]{\mbox{$G^{#1}_{#2}(#3,#4)$}}
%%%%%%%%%%%
\newtheorem{theorem}{Theorem}[section]
\newtheorem{corollary}[theorem]{Corollary}
\newtheorem{lemma}[theorem]{Lemma}
\newtheorem{definition}[theorem]{Definition}
\newtheorem{remark}[theorem]{Remark}
\newtheorem{remarks}[theorem]{Remarks}
\newtheorem{proposition}[theorem]{Proposition}
\newtheorem{problem}{Problem}
%%%%%%%%%%%%%%%%%%
\newtheorem{conjecture}[theorem]{Conjecture}
\newtheorem{question}[theorem]{Question}
\newtheorem{example}[theorem]{Example}
%%%%%%%%%%%%%%%%%%%
\newtheorem{Thm}[theorem]{Theorem}
\newtheorem{Lem}[theorem]{Lemma}
\newtheorem{Pro}[theorem]{Proposition}
\newtheorem{Def}[theorem]{Definition}
\newtheorem{Exa}[theorem]{Example}
\newtheorem{Exs}[theorem]{Examples}
\newtheorem{Rems}[theorem]{Remarks}
\newtheorem{rem}[theorem]{Remark}
\newtheorem{Cor}[theorem]{Corollary}
\newtheorem{Conj}[theorem]{Conjecture}
\newtheorem{Prob}[theorem]{Problem}
\newtheorem{Ques}[theorem]{Question}
%%%%%%%%%%%%%
\newcommand{\pf}{\noindent \mbox{{\bf Proof}: }}
\newcommand{\dnorm}[1]{\thkl #1 \thkl\,}
%%%%%%%%%%%%%%%%%%
%\newenvironment{proof}{{\bf Proof.}}{\hfill $\bowtie$\vskip4mm}

\renewcommand{\theequation}{\thesection.\arabic{equation}}
\catcode`@=11 \@addtoreset{equation}{section} \catcode`@=12
%%%%%%%%%%%%%%%%%
%%Macros for Greek letters.
\def\ga{\alpha}     \def\gb{\beta}       \def\gg{\gamma}
\def\gc{\chi}       \def\gd{\delta}      \def\ge{\epsilon}
\def\gth{\theta}                         \def\vge{\varepsilon}
       \def\vgf{\varphi}    \def\gh{\eta}
\def\gi{\iota}      \def\gk{\kappa}      \def\gl{\lambda}
\def\gm{\mu}        \def\gn{\nu}         \def\gp{\pi}
\def\vgp{\varpi}    \def\gr{\rho}        \def\vgr{\varrho}
\def\gs{\sigma}     \def\vgs{\varsigma}  \def\gt{\tau}
\def\gu{\upsilon}   \def\gv{\vartheta}   \def\gw{\omega}
        \def\gy{\psi}        \def\gz{\zeta}
\def\Gg{\Gamma}     \def\Gd{\Delta}      \def\Gf{\Phi}
\def\Gth{Theta}
\def\Gl{\Lambda}    \def\Gs{\Sigma}      \def\Gp{\Pi}
\def\Gw{\Omega}     \def\Gx{\Xi}         \def\Gy{\Psi}
%%%%%%%%%%%%%%%%%%%%%%%%%%%%%%%%

%\begin{titlepage}

\title{A Liouville-type theorem for Schr\"odinger operators}
\author{Yehuda Pinchover\\
 {\small Department of Mathematics}\\ {\small  Technion - Israel Institute of Technology}\\
 {\small Haifa 32000, Israel}\\
{\small pincho@techunix.technion.ac.il}}
\date{}
\maketitle
\begin{abstract} In this paper we prove a sufficient condition, in terms of
the behavior of a ground state of a symmetric critical operator
$P_1$, such that a nonzero subsolution of a symmetric nonnegative
operator $P_0$ is a ground state. Particularly, if
$P_j:=-\Delta+V_j$, for $j=0,1$, are two nonnegative Schr\"odinger
operators defined on $\Omega\subseteq \mathbb{R}^d$ such that
$P_1$ is critical in $\Omega$ with a ground state $\varphi$, the
function $\psi\nleq 0$ is a subsolution of the equation $P_0u=0$
in $\Omega$ and satisfies $|\psi|\leq C\varphi$ in $\Omega$, then
$P_0$ is critical in $\Omega$ and $\psi$ is its ground state. In
particular, $\psi$ is (up to a multiplicative constant) the unique
positive supersolution of the equation $P_0u=0$ in $\Omega$.
Similar results hold for general symmetric operators, and also on
Riemannian manifolds.
\\[1mm]
\end{abstract}
%%%%%%%%%%%%%%%%%%%%%%%%%%%%%%%%%%
%\end{titlepage}
%%%%%%%%%%%%%%%%%%%%%%%%%%%%%%%%%%%%%%%%%%%%%%%%%
%\tableofcontents
%\newpage
%%%%%%%%%%%%%%%%%%%%%%%%%%%%%%%%%%
 \mysection{Introduction}\label{secint}
%%%%%%%%%%%%%%%%%%%%%%%%%%%%%%%%%%%%%%%%%%%%%%
Let $\Omega\subset\R^d$ be a domain.  We assume that $A:\Omega
\rightarrow \mathbb{R}^{d^2}$ is a measurable matrix valued
function such that for every compact set $K\subset \Omega$  there
exists $\mu_K>1$ such that \be \label{stell} \mu_K^{-1}I_d\le
A(x)\le \mu_K I_d \qquad \forall x\in K,
 \ee
where $I_d$ is the $d$-dimensional identity matrix, and the matrix
inequality $A\leq B$ means that $B-A$ is a nonnegative matrix on
$\mathbb{R}^d$. Let $V\in
L^{p}_{\mathrm{loc}}(\Omega;\mathbb{R})$, where $p>{d}/{2}$. We
consider the quadratic form \be \label{assume}
\mathbf{a}_{A,V}[u]:=\int_\Omega\left(A\nabla u\cdot \nabla
u+V|u|^2\right)\mathrm{d}x   \ee on $\core$
 associated with the Schr\"odinger
equation \be \label{divform}
Pu:=(-\nabla\cdot(A\nabla)+V)u=0\qquad \mbox{ in } \Omega. \ee
 We say that $\mathbf{a}_{A,V}$ is {\em nonnegative} on $C_0^\infty(\Omega)$,
 if $\mathbf{a}_{A,V}[u]\geq 0$ for all $u\in \core$.
\begin{definition}{\em We say that $v\in H^{1}_{\mathrm{loc}}(\Omega)$ is a {\em (weak) solution} of
\eqref{divform}  if for every $\varphi\in\core$
 \be \label{solution} \int_\Omega (A\nabla v\cdot\nabla
\varphi+Vv\varphi)\dx=0. \ee We say that $v\in
H^{1}_{\mathrm{loc}}(\Omega)$ is a {\em subsolution} of
\eqref{divform} if for every nonnegative $\varphi\in\core$
 \be\label{subsolution}
\int_\Omega (A\nabla v\cdot\nabla \varphi+Vv\varphi)\dx\leq 0.
 \ee
$v\in
H^{1}_{\mathrm{loc}}(\Omega)$ is a {\em supersolution} of
 \eqref{divform} if $-v$ is a subsolution of \eqref{divform}.
}\end{definition}
 Let
$\mathcal{C}_P(\Omega)$ be the cone of all positive solutions of
the equation $Pu=0$ in $\Omega$, and let
\begin{equation}\label{lambda0}
\lambda_0(P,\Omega) := \sup\{\lambda \in \mathbb{R} \mid
\mathcal{C}_{P-\lambda}(\Omega)\neq \emptyset\}
\end{equation}
be the {\em generalized principal eigenvalue} of the operator $P$
in $\Omega$. By the Allegretto-Piepenbrink theory (see for
example, \cite{Agmon82,Pins95}), the form $\mathbf{a}_{A,V}$ is
nonnegative on $\core$ if and only if $\lambda_0(P,\Omega)\geq 0$.

Let $K\Subset \Omega$ (i.e. $K$ is relatively compact in
$\Omega$). Recall \cite{Agmon82,Pins95} that $u\in
\mathcal{C}_P(\Omega\setminus K)$ is said to be a {\em positive
solution of the operator $P$ of minimal growth in a neighborhood
of infinity in} $\Omega$, if for  any $K\Subset K_1 \Subset
\Omega$, with a smooth boundary,  and any $v\in
\mathcal{C}_P(\Omega\setminus K_1)\cap C((\Omega\setminus K_1)\cup
\partial K_1)$, the inequality $u\le v$
on $\partial K_1$ implies that $u\le v$ in $\Omega\setminus K_1$.
A positive solution $u\in \mathcal{C}_P(\Omega)$ which has minimal
growth in a neighborhood of infinity in $\Omega$ is called a {\em
ground state} of $P$ in $\Omega$.

The operator $P$ is said to be {\em critical} in $\Omega$, if $P$
admits a ground state in $\Omega$. The operator $P$ is called {\em
subcritical} in $\Omega$, if $\mathcal{C}_P(\Omega)\neq
\emptyset$, but $P$ is not critical in $\Omega$. If
$\mathcal{C}_P(\Omega)= \emptyset$, then $P$ is {\em
supercritical} in $\Omega$.

It is known that the operator $P$ is critical in $\Omega$ if and
only if the equation $Pu=0$ in $\Omega$ admits (up to a
multiplicative constant) a unique positive supersolution. In
particular, in the critical case we have $\dim
\mathcal{C}_{P}(\Omega)=1$ (see for example \cite{P2006,Pins95}
and the references therein).

On the other hand, $P$ is subcritical in $\Omega$, if and only if
$P$ admits a {\em positive minimal Green function}
$G^\Omega_P(x,y)$ in $\Omega$. For each $y\in \Omega$, the
function $G^\Omega_P(\cdot,y)$ is a positive solution of the
equation $Pu=0$ in $\Omega\setminus\{y\}$ that has minimal growth
in a neighborhood of infinity in $\Omega$ and has a (suitably
normalized) nonremovable singularity at $y$ (see for example
\cite{P2006,Pins95} and the references therein).

The following basic example will be used few times along the
paper.
\begin{example}\label{ex0} {\em Let $P\!=\!-\Delta$ and $\Omega\!=\!\mathbb{R}^d$. It is well known that
$\lambda_0(-\Delta,\mathbb{R}^d)\!=\!0$. In addition, the positive
Liouville theorem asserts that
$$\mathcal{C}_{-\Delta}(\mathbb{R}^d)=\{c\mathbf{1}\mid c>0\},$$
where $\mathbf{1}$ is the constant function taking at any point
the value $1$. Moreover, $-\Delta$ is critical in $\mathbb{R}^d$
if and only if $d\leq 2$.
 }\end{example}
 Recently, Berestycki, Hamel, and Roques \cite{BHR} has
introduced the following definition which arises  naturally in the
study of some semilinear equations.
\begin{definition}{\em
\begin{multline}
\lambda_0'(P,\Omega) := \inf\{\lambda \in \mathbb{R} \mid \exists
\phi\in C^2(\Omega)\cap W^{2,\infty}(\Omega), \phi>0,
(P-\lambda)\phi\leq 0 \mbox{ in } \Omega,\nonumber \\\phi=0 \mbox{
on } \partial \Omega, \mbox{ if }
\partial \Omega\neq \emptyset \}.\label{lambda0'}
\end{multline}
 }\end{definition}

\vskip 2mm

Before formulating our results, we present four basic problems
(Problems~\ref{problem1}--\ref{problem4}) concerning the
particular case $\Omega= \mathbb{R}^d$, which have been solved in
the past few years. It turns out that (almost) all these results
follow directly from our main result (Theorem~\ref{mainthm}).

\begin{problem}[\cite{S81,S82}]\label{problem1}
Let $V\in L^2_{\mathrm{loc}}(\mathbb{R}^d)$. Does the existence of
a positive bounded solution to the equation
 \be\label{eqhq}
H_Vu:=(-\Delta+V)u=0\qquad  \mbox {on } \mathbb{R}^d
 \ee
imply that $H_V$ is critical in $\mathbb{R}^d$?
\end{problem}
\begin{problem}[\cite{BCN}]\label{problem2}
Suppose that $V$ is smooth and bounded. Does the existence of a
sign-changing bounded solution to equation \eqref{eqhq} imply that
$\lambda_0(H_V,\mathbb{R}^d)<0$?
\end{problem}
\begin{problem}[\cite{BCN,GG}]\label{problem3}
Let $\sigma$ be a strictly positive $C^2$-function on
$\mathbb{R}^d$, and consider the divergence form operator $L=
\nabla\cdot (\sigma^2\nabla)$  on $\mathbb{R}^d$. Suppose that the
equation $L\psi=0$ in $\mathbb{R}^d$ admits a nonzero solution
$\psi$ such that $\psi\sigma$ is bounded.  Is $\psi$ necessarily
the constant function?
\end{problem}
\begin{problem}[{\cite[Conjecture~4.6]{BR}}]\label{problem4}
Suppose that $P=-\nabla\cdot(A\nabla)+V$ is uniformly elliptic
operator with smooth bounded coefficients on $\mathbb{R}^d$. Does
the inequality
 $$
\lambda_0(P,\mathbb{R}^d)\leq \lambda_0'(P,\mathbb{R}^d)
 $$
holds true in any dimension $d$.
\end{problem}

\vskip 2mm

The answers to the above four problems for the free Laplacian in
$\mathbb{R}^d$ are well known. Nevertheless, the above problems
are not of perturbational nature since there is no assumption on
the asymptotic behavior of the coefficients of the given operator
near infinity.

Problem~\ref{problem1} was posed by B.~Simon in \cite{S81,S82}.
Clearly, the answer  to Problem~\ref{problem1} is `no' for $d\geq
3$. Partial results concerning Problem~\ref{problem1} for $d\leq
2$ were obtained under the assumption that $V$ is a short-range
potential (see for example, \cite{GZ1,GZ2,M84,P2006}). On the
other hand, Gesztesy, and Zhao showed in \cite[Example~4.6]{GZ1}
that there is a critical Schr\"odinger operator on $\mathbb{R}$
with `almost' short-range potential such that its ground state
behaves logarithmically.

In a recent article Damanik, Killip, and Simon proved a result
which reveals a complete answer to Problem~\ref{problem1}.
\begin{theorem}[cf. {\cite[Theorem~5]{DKS}}]\label{thmDKS}
The answer to Problem~\ref{problem1} is ``yes" if and only if
$d=1,2$.
\end{theorem}

Indeed, for $d = 1,2$, it is shown in \cite{DKS} that  if the
equation $H_Vu=0$ admits a positive bounded solution, then any
$W\in L^2_{\mathrm{loc}}(\mathbb{R}^d)$ satisfying $H_{V\pm W}\geq
0$ is necessarily the zero potential. But this property holds if
and only if $H_V$ is critical (see \cite{P2006}).

\vskip 4mm

Let us turn to Problem~\ref{problem2} which was raised by
Berestycki, Caffarelli, and Nirenberg \cite{BCN}. This problem is
closely related to De Giorgi's conjecture \cite{DG} (see
\cite{B,BBG,BCN,GG}). In \cite{GG}, Ghoussoub and Gui showed a
connection between Problem~\ref{problem2} and
Problem~\ref{problem3} which concerns the Liouville property for
operators in divergence form (see also the proof of Theorem 1.7 in
\cite{BCN}). In fact, the following result is proved in
\cite{BCN,GG,B}.
\begin{theorem}\label{thmBCNGGB}
The answers to problems~\ref{problem2} and \ref{problem3} are
``yes" if and only if $d=1,2$.
\end{theorem}
Note that Ghoussoub and Gui \cite{GG} used this Liouville-type
theorem for $d=2$ \cite{BCN}, to prove De Giorgi's Conjecture in
$\mathbb{R}^2$.

\vskip 4mm

Problem~\ref{problem3} was posed by Berestycki and Rossi \cite{BR}
who also
 solved it for $d\leq 3$:
\begin{theorem}[{\cite[Theorem~4.1]{BR}}]\label{thmBR}
Suppose that $P=-\nabla\cdot(A\nabla)+V$ is uniformly elliptic
operator with smooth bounded coefficients on $\mathbb{R}^d$.  If
$d\leq 3$, then
 $$
\lambda_0(P,\mathbb{R}^d)\leq \lambda_0'(P,\mathbb{R}^d).
 $$
\end{theorem}

\vskip 2mm

The purpose of the present article is to (partially) generalize
theorems~\ref{thmDKS}, \ref{thmBCNGGB}, and \ref{thmBR} to general
symmetric operators which are defined on an arbitrary domain
$\Omega \subseteq \mathbb{R}^d$, or on a noncompact Riemannian
manifold. Our main result is as follows.
\begin{theorem}\label{mainthm} Let $\Omega$ be a domain in $\R^d$, $d\geq 1$.
Consider two Schr\"odinger operators defined on $\Omega$ of the
form
\begin{equation}\label{eqpj}
P_j:=-\nabla\cdot(A_j\nabla)+V_j\qquad j=0,1,
\end{equation}
such that  $V_j\in L^{p}_{\mathrm{loc}}(\Omega;\mathbb{R})$ for
some $p>{d}/{2}$, and $A_j$ satisfy \eqref{stell}.

 Assume that the following assumptions hold true.

\begin{itemize}
\item[(i)] The operator  $P_1$ is critical in $\Omega$. Denote
by $\varphi\in \mathcal{C}_{P_1}(\Omega)$ its ground state.

\item[(ii)]  $\lambda_0(P_0,\Omega)\geq 0$, and there exists a
real function $\psi\in H^1_{\mathrm{loc}}(\Omega)$ such that
$\psi_+\neq 0$, and $P_0\psi \leq 0$ in $\Omega$, where
$u_+(x):=\max\{0, u(x)\}$.

\item[(iii)] The following matrix inequality holds
\begin{equation}\label{psialephia}
(\psi_+)^2(x) A_0(x)\leq C\varphi^2(x) A_1(x)\qquad  \mbox{ a. e.
in } \Omega,
\end{equation}
where $C>0$ is a positive constant.
\end{itemize}
Then the operator $P_0$ is critical in $\Omega$, and $\psi$ is its
ground state. In particular, $\dim \mathcal{C}_{P_0}(\Omega)=1$
and $\lambda_0(P_0,\Omega)=0$.
\end{theorem}
\begin{corollary}\label{corbcn}
Suppose that all the assumptions of Theorem~\ref{mainthm} are
satisfied except possibly the assumption that
$\lambda_0(P_0,\Omega)\geq 0$. Assume further that either $\psi$
changes its sign in $\Omega$, or $\psi$ is not a solution  of the
equation $P_0u=0$ in $\Omega$. Then $\lambda_0(P_0,\Omega)<0$.
\end{corollary}

Theorem~\ref{mainthm} and Corollary~\ref{corbcn} imply in
particular the sufficiency parts of Theorem~\ref{thmDKS} and
Theorem~\ref{thmBCNGGB}, and also Theorem~\ref{thmBR} for $d=1,2$.
Note that in contrast to the assumptions of
theorem~\ref{thmBCNGGB} and \ref{thmBR}, we assume in
Theorem~\ref{mainthm} neither that the functions $V_j$ are bounded
and smooth, nor that the matrix valued functions $A_j$ are smooth.

The outline of the paper is as follows. In
Section~\ref{secpreliminaries}, we present some results from
\cite{PT2} that will be used in the proof of
Theorem~\ref{mainthm}. Section~\ref{secpf} is devoted to the proof
of Theorem~\ref{mainthm} and its consequences. We conclude the
paper in Section~\ref{secopen}, where we pose two open problems
suggested by the results of the present paper.
%%%%%%%%%%%%%%%%%%%%%%%%%%%%%%%%%%%%
\mysection{Preliminary results}\label{secpreliminaries}
%%%%%%%%%%%%%%%%%%%%%%%%%%%%%%%%%%%%%%%%%%%%%%
\begin{definition}\label{defnull}{\em We say that a sequence
$\{u_k\}\subset\core$ is a {\em null sequence} for the nonnegative
quadratic form $\mathbf{a}_{A,V}$ if there exists an open set
$B\Subset\Omega$ such that $\int_B|u_k|^2\dx=1$, and
 \be
\lim_{k\to\infty}\mathbf{a}_{A,V}[u_k]=0.\ee We say that a
positive function $\varphi$ is a {\em null state} for the
nonnegative quadratic form $\mathbf{a}_{A,V}$, if there exists a
null sequence $\{u_k\}$ for the form $\mathbf{a}_{A,V}$ such that
$u_k\to\varphi$ in $L^2_\mathrm{loc}(\Omega)$.
 }\end{definition}
\begin{remark}\label{remc1}{\em
The requirement that $u_k\subset \core$, can clearly be weakened
by assuming only that $\{u_k\}\subset H^1_0(\Omega)$. Also, the
requirement that $\int_B|u_k|^2\dx=1$ can be replaced by
$\int_B|u_k|^2\dx\asymp 1$, where $f_k\asymp g_k$ means that there
exists a positive constant $C$ such that $C^{-1}g_k\leq f_k \leq
Cg_k$  for all $k\in \mathbb{N}$.
 }\end{remark}

%%%%%%%%%%%%%%%%%%%%%%%%%%%%%%%%%%%%%
The following auxiliary lemma is well known (see, e.g.
\cite{DS84,M2002,PT2}).
\begin{lemma}\label{lemqf}  Let $\psi\in H^1_{\mathrm{loc}}(\Omega)$ be a nonnegative subsolution of the
equation $P\psi=0$ in $\Omega$. Then for any nonnegative $v\in
\core$ we have
 \begin{equation}
\label{e1}\mathbf{a}_{A,V}[\psi v]\leq
 \int_\Omega (\psi)^2
A\nabla v\cdot\nabla v \,\mathrm{d}x.
 \end{equation}
Moreover, if $\psi$ is a (real valued) solution of the equation
$P\psi=0$ in $\Omega$, then for any $v\in \core$ we have
\begin{equation}
\label{e2}\mathbf{a}_{A,V}[\psi v]=
 \int_\Omega (\psi)^2
A\nabla v\cdot\nabla v \,\mathrm{d}x.
 \end{equation}
  \end{lemma}
\begin{proof} Follows from the definition of a weak (sub)solution and elementary calculation.
\end{proof}
The following theorem was recently proved by K.~Tintarev and the
author \cite{PT2} (see also \cite{PT3}).
\begin{theorem}\label{thmnull}
Suppose that $\mathbf{a}_{A,V}\geq 0$ on $\core$. Then
$\mathbf{a}_{A,V}$ has a null sequence if and only if the
corresponding operator $P=-\nabla\cdot(A\nabla)+V$ is critical in
$\Omega$. In this case, any null sequence converges in
$L^2_{\mathrm{loc}}(\Omega)$ to $c\varphi$, where $\varphi$ is a
ground state of the operator $P$ and $c$ is a nonzero constant.

Moreover, there exists a null sequence $\{u_k\}$ of nonnegative
functions that converges to $\varphi$ locally uniformly in
$\Omega\setminus\{x_0\}$, where $x_0$ is some point in $\Omega$.
\end{theorem}
\section{Proof of Theorem~\ref{mainthm}}\label{secpf}
%%%%%%%%%
In this section we prove Theorem~\ref{mainthm} and some
consequences.
\begin{proof}[Proof of Theorem~\ref{mainthm}]
Since $\psi$ satisfy $P_0\psi\leq 0$ in $\Omega$, it follows that
$P_0\psi_+\leq 0$ in $\Omega$ (see for example
\cite[Lemma~2.7]{Agmon82}).

By Theorem~\ref{thmnull} and our assumptions, there exists a null
sequence $\{u_k\}$ for the quadratic form $\mathbf{a}_{A_1,V_1}$
of nonnegative functions  which converges locally uniformly in
$\Omega\setminus\{x_0\}$ and in $L^2_{\mathrm{loc}}(\Omega)$ to
the ground state $\varphi$ of the operator $P_1$, and satisfies
$\int_B(u_k)^2\dx=1$ for some open set
$B\Subset\Omega\setminus\{x_0\}$ and all $k\in \mathbb{N}$.

Denote $w_k:=u_k/\varphi$. Since $w_k\to \mathrm{constant}$
locally uniformly in $\Omega\setminus\{x_0\}$ and $\psi_+\neq 0$,
it follows that $\int_{B_1}(\psi_+ w_k)^2\dx\asymp 1$ for some
open set $B_1\Subset\Omega$ and every $k\geq k_0$. Moreover, by
Lemma~\ref{lemqf} and our assumptions, we have
\begin{multline}\label{eq1}
\mathbf{a}_{A_0,V_0}[\psi_+ w_k]\leq \int_\Omega
(\psi_+)^2A_0\nabla w_k\cdot\nabla w_k
\,\mathrm{d}x\leq\\[4mm]
C\int_\Omega \varphi^2A_1\nabla w_k\cdot\nabla w_k \,\mathrm{d}x=
C\mathbf{a}_{A_1,V_1}[\varphi w_k]=C\mathbf{a}_{A_1,V_1}[u_k]\to
0.
\end{multline}
 Therefore,
$\{\psi_+ w_k\}$ is a null sequence for $P_0$. By
Theorem~\ref{thmnull}, $P_0$ is critical in $\Omega$ and $\psi_+$
is its ground state. In particular, $\psi_+$ is strictly positive,
and hence $\psi_-=0$, and $\psi=\psi_+$ is the ground state of
$P_0$.
\end{proof}
%%%%
\begin{remark}\label{rem22} {\em
Suppose that all the assumptions of Theorem~\ref{mainthm} are
satisfied except possibly the assumption that
$\lambda_0(P_0,\Omega)\geq 0$. One can show directly that
$\lambda_0(P_0,\Omega)\leq 0$. Indeed, using the notations of the
proof of Theorem~\ref{mainthm}, we have that for some $C_1>0$
$$\int_\Omega(\psi_+w_k)^2\dx\geq C_1\int_B(u_k)^2\dx=C_1\qquad \forall k\geq k_0.$$
Moreover, by Lemma~\ref{lemqf} and our assumptions, we have
\begin{multline}\label{eq7}
\dfrac{\mathbf{a}_{A_0,V_0}[\psi_+
u_k]}{\int_\Omega(\psi_+w_k)^2\dx}\leq \dfrac{\int_\Omega
(\psi_+)^2A_0\nabla w_k\cdot\nabla w_k
\,\mathrm{d}x}{\int_\Omega(\psi_+w_k)^2\dx}\leq\\[4mm]
\tilde{C}\dfrac{\int_\Omega \varphi^2A_1\nabla w_k\cdot\nabla w_k
\,\mathrm{d}x}{\int_B(u_k)^2\dx}=
\tilde{C}\mathbf{a}_{A_1,V_1}[\varphi
w_k]=\tilde{C}\mathbf{a}_{A_1,V_1}[u_k]\to 0.
\end{multline}
Therefore, the Rayleigh-Ritz variational formula implies that
$\lambda_0(P_0,\Omega)\leq 0$.

It follows that  \bean \lambda_0(P_0,\Omega)\leq \inf\{\lambda \in
\mathbb{R} &\mid&
\exists\, \psi\nleqq 0, (P_0-\lambda)\psi\leq 0 \mbox{ in } \Omega \mbox{ s.t. } \\
 &&\psi^2(x) A_0(x)\leq C\varphi^2(x) A_1(x) \mbox{ in }
\Omega \mbox{ for some }
\\   &&\mbox{critical operator $P_1$ with a ground state } \varphi
\}.
\eean
 In particular, if $P=-\nabla\cdot(A\nabla)+V$ is an elliptic operator on  $\mathbb{R}^d$, $d\leq
2$, with a bounded matrix $A$, then $\lambda_0(P,\mathbb{R}^d)\leq
\lambda_0'(P,\mathbb{R}^d)$ (cf. Theorem~\ref{thmBR}).
 }\end{remark}
Recall that if $P:=-\nabla\cdot(A\nabla)+V$ is
$\mathbb{Z}^d$-periodic on $\mathbb{R}^d$, then  $P-\lambda_0$
admits a (unique) periodic positive solution (see for example
\cite{LP,Pins95}). On the other hand,  $-\Delta$ is critical in
$\mathbb{R}^d$ if and only if $d\leq 2$ (see Example~\ref{ex0}).
Therefore, Theorem~\ref{mainthm} implies the following result of
R.~Pinsky (who proved it for general second-order elliptic
$\mathbb{Z}^d$-periodic operators).
\begin{corollary}[\cite{Pinsper}]\label{corper}
Assume that the coefficients of the elliptic operator
$P:=-\nabla\cdot(A\nabla)+V$ are $\mathbb{Z}^d$-periodic on
$\mathbb{R}^d$. Then the operator $P-\lambda_0$ is critical in
$\mathbb{R}^d$ if and only if $d\leq 2$.
\end{corollary}
\begin{remark}\label{rem2} {\em
Suppose that $P_j$ are two nonnegative symmetric operators which
are defined on a noncompact Riemannian manifold $M$ of dimension
$d$, where $j=0,1$. Since Lemma~\ref{lemqf} holds true also in
this case (see \cite{M2002}), it follows that
Theorem~\ref{thmnull} is valid on Riemannian manifolds, which in
turn implies that Theorem~\ref{mainthm} holds true also in this
case.
 }\end{remark}
Recall that a Riemannian manifold $M$ is called {\em recurrent} if
the Laplace-Beltrami operator on $M$ is critical (see
\cite{Pins95}). Therefore, we have in particular, the following
generalization of Theorem~\ref{thmDKS} and
Theorem~\ref{thmBCNGGB}.
\begin{theorem}\label{thmDKSR}
Let $M$ be a recurrent Riemannian noncompact manifold  of
dimension $d$. Let $V \in L^2_{\mathrm{loc}}(M)$. Suppose that
$H_V:=-\Delta+V \geq 0$ on $C_0^\infty(M)$, and that the equation
$H_Vu=0$ in $M$ admits a nonzero bounded subsolution $\psi$ such
that $\psi_+\neq 0$. Then $H_V$ is critical in $M$ and $\psi$ is a
ground state of $H_V$ in $M$. In particular, $\lambda_0(H_V)=0$,
the space of all bounded solutions of the equation $H_Vu=0$ in $M$
is one-dimensional, and $\dim \mathcal{C}_{H_V}(M)=1$ .
\end{theorem}
In addition, one can use the results in \cite{LP} and
\cite[Theorem~5.2.11]{D89} to extend Corollary~\ref{corper} to the
case of equivariant Schr\"odinger operators on cocompact nilpotent
coverings.
\begin{corollary}\label{cornilp}
 Let $M$ be a noncompact
nilpotent covering of a compact Riemannian manifold of dimension
$d$. Suppose that $P:=-\Delta+V$ is an equivariant operator on $M$
with respect to the (nilpotent) deck group $G$. Then $P-\lambda_0$
is critical in $M$ if and only if $G$ has a normal subgroup of
finite index isomorphic to $\mathbb{Z}^n$ for $n\leq 2$.
 \end{corollary}
 Theorem~\ref{mainthm} can be considered as a sufficient
condition for the removability of singularity at infinity in
$\Omega$ or as a Phragm\'{e}n-Lindel\"of type principle. A
positive solution of \eqref{divform} in $\Omega\setminus K$, where
$K\Subset\Omega$, is called {\em singular at infinity} if it does
not have minimal growth in a neighborhood of infinity in $\Omega$.
Accordingly, Theorem~\ref{mainthm} implies that the behavior of a
positive solution of minimal growth in a neighborhood of infinity
in $\Omega$ of an equation of the form \eqref{eqpj}, dictates some
`growth' on all positive singular at infinity solutions of {\em
any} equation of the form \eqref{eqpj}. More precisely, we have
the following result.
\begin{corollary}\label{cormin}
Suppose that for $j=0,1$ the operators $P_j$ are of the form
\eqref{eqpj}, and $A_j$ satisfy \eqref{stell}. Let $u_1$ be a
positive solution of the equation $P_1u=0$ of minimal growth in a
neighborhood of infinity in $\Omega$, and let $u_0$ be a positive
solution of the equation $P_0u=0$ in $\Omega\setminus K$, where
$K\Subset\Omega$. If $(u_0)^2A_0\leq C (u_1)^2A_1$ in
$\Omega\setminus K$, then $u_0$ is nonsingular at infinity, i.e.,
$u_0$ is a positive solution of the equation $P_0u=0$ of minimal
growth in a neighborhood of infinity in $\Omega$.
\end{corollary}

\begin{proof}[Proof of Corollary~\ref{cormin}] Let
$\widetilde{u_0},\widetilde{u_1}\in H^{1}_{\mathrm{loc}}(\Omega)$
be positive functions which are defined in $\Omega$ such that
$\widetilde{u_j}|_{\Omega\setminus K_1}=u_j$, and
$\widetilde{u_j}|_{\overline{K_1}}$ are sufficiently smooth, where
$K_1\Subset\Omega$, and $j=0,1$.

Then for $j=0,1$, $\widetilde{u_j}\in
\mathcal{C}_{\widetilde{P_j}}(\Omega)$, where  the operators
$\widetilde{P_j}$ are of the form \eqref{eqpj} and satisfy
$\widetilde{P_j}|_{\Omega\setminus K_2}=P_j$ for some $K_2\Subset
\Omega$. Since $u_1$ (and hence also $\widetilde{u_1}$) is a
positive solution of the equation $\widetilde{P_1}u=0$ of minimal
growth in a neighborhood of infinity in $\Omega$, it follows that
$\widetilde{u_1}$ is a ground state of the critical operator
$\widetilde{P_1}$  in $\Omega$. Therefore, Theorem~\ref{mainthm}
implies that $\widetilde{u_0}$ is a ground state of the critical
operator $\widetilde{P_0}$ in $\Omega$. Hence, $u_0$ is a positive
solution of the equation $P_0u=0$ of minimal growth in a
neighborhood of infinity in $\Omega$.
\end{proof}
\begin{example}\label{ex1} {\em
Let $d\geq 2$, and $V\in L^p_{\mathrm{loc}}(\mathbb{R}^d)$, where
$p>{d}/{2}$.  Suppose that $H_V:=-\Delta+V\geq 0$ on
$C_0^\infty(\mathbb{R}^d)$, and the equation $H_Vu=0$ on
$\mathbb{R}^d$ has a subsolution solution $\psi\nleq 0$ satisfying
\be \psi_+(x)=O(|x|^{\frac{2-d}{2}})\qquad \mbox{ as } |x| \to
\infty.
 \ee
Since $\varphi(x):= |x|^{\frac{2-d}{2}}$ is a positive solution of
the Hardy-type equation
$$-\Delta u-\left(\frac{d-2}{2}\right)^2\frac{u}{|x|^{2}}=0$$
of minimal growth in a neighborhood of infinity in $\mathbb{R}^d$,
it follows from Corollary~\ref{cormin} that $H_V$ is critical in
$\mathbb{R}^d$ and $\psi$ is its ground state (cf. Theorem~1.7 in
\cite{BCN}).
 }\end{example}
\begin{example}\label{ex7} {\em
Let $d=1$, and $V\in L^p_{\mathrm{loc}}(\mathbb{R})$, where $p>1$.
Suppose that $H_V:=-\mathrm{d}^2/\mathrm{d}x^2+V\geq 0$ on
$C_0^\infty(\mathbb{R})$, and the equation $H_Vu=0$ on
$\mathbb{R}$ has a subsolution solution $\psi\nleq 0$ satisfying
\be \psi_+(x)=O(\log |x|)\qquad \mbox{ as } |x| \to \infty.
 \ee
It follows from \cite[Example~4.6]{GZ1} and Corollary~\ref{cormin}
that $H_V$ is critical in $\mathbb{R}$  and $\psi$ is its ground
state.
 }\end{example}
%%%%%%%%%%%%%%%%%%%%%%%%%%%%%%
\mysection{Open problems}\label{secopen} We conclude the paper
with two open problems suggested by the above results which are
left for future investigation.
\begin{problem}\label{pr1}
Generalize Theorem~\ref{mainthm} to the class of nonsymmetric
second-order linear elliptic operators with real coefficients
which have the same (or even comparable) principal parts, or at
least to the subclass of operators which differ only by their
zero-order terms.
\end{problem}

\begin{remarks}\label{rem1} {\em
1. Clearly, the condition \eqref{psialephia} which involves the
{\em squares} of the corresponding solutions of the symmetric
operators $P_j$, for $j=0,1$, should be replaced in the
nonsymmetric case by a condition which involves {\em products} of
the form $u_ju_j^*\,$, where $u_j$ (resp. $u_j^*$) are the
corresponding solutions of the operators $P_j$ (resp. of the
formal adjoint operators $P^*_J$) for $j=0,1$.

2. Let $u$ be a positive solution of an equation of the form
\eqref{divform} of minimal growth in a neighborhood of infinity in
$\Omega$, then Corollary~\ref{cormin} implies that any positive
solution $v$ of another equation of the form \eqref{divform} (with
a comparable principal part) in a neighborhood of infinity in
$\Omega$ which is singular at infinity satisfies

$$\liminf_{x\to \infty} \frac{u(x)}{v(x)}=0$$
in the one-point compactification of $\Omega$ ($\infty$ denotes
the point at infinity in $\Omega$).

Ancona proved \cite{An} that a subcritical symmetric (or even
quasi-symmetric) operator $P$ in $\Omega$ always admits $v\in
\mathcal{C}_{P}(\Omega)$, such that
$$\lim_{x\to \infty} \frac{G_P^\Omega(x, x_0)}{v(x)}=0.$$
Moreover, it is shown in \cite{An}  that such a positive solution
does not always exist for general nonsymmetric operators. This
result indicates that the answer to Problem~\ref{pr1} in the
nonsymmetric case might be more involved.
 }\end{remarks}
The second problem that we pose deals with Liouville-type theorems
for the $p$-Laplacian with a potential term. Let $\Omega$ be a
domain in $\mathbb{R}^d$, $d\geq 2$, and $1<p<\infty$. Fix $V\in
L_{\mathrm{loc}}^\infty(\Omega)$. Recently the criticality theory
for linear equations was extended in \cite{PT3} to quasilinear
equations of the form
 \be\label{plaplace}
 -\nabla\cdot(|\nabla u|^{p-2}\nabla u)+V|u|^{p-2}u=0\qquad \mbox{ in }
 \Omega.
 \ee
In particular, Theorem~\ref{thmnull} was proved also for such
equations. Therefore, it is natural to pose the following problem.
\begin{problem}\label{pr2}
Assume that $1<p \leq d$. Generalize Theorem~\ref{mainthm} to
positive solutions of quasilinear equations of the form
\eqref{plaplace}.
\end{problem}
%%%%%%%%%%%%%%%%%%%%%%%%%%%%%%%%%%%%%%%%%%%%%%%%%%%%%%%%%%%%%%%%%%
\begin{center}
{\bf Acknowledgments} \end{center} The author wishes to thank
H.~Brezis, F.~Gesztesy, and M.~Marcus for valuable discussions.
This work was partially supported by the RTN network ``Nonlinear
Partial Differential Equations Describing Front Propagation and
Other Singular Phenomena", HPRN-CT-2002-00274, the Israel Science
Foundation founded by the Israeli Academy of Sciences and
Humanities, and the Fund for the Promotion of Research at the
Technion.
%%%%%%%%%%%%%%%%%%%%%%%%%%%%%%%%%%%%%%%%%%%%%%%%%

\end{document}